\documentclass{article}

\usepackage[tbtags]{amsmath}
\usepackage{amsfonts,amssymb}
\usepackage{theorem}

\def\inv{\mathrm{inv}}
\def\const{\operatorname{const}}
\newcommand\C{\mathbb{C}}

\begin{document}

\title{Higher discriminants of binary forms\footnote{Published in
{\it Theor.~Math.~Phys.}~(Nov 2007) {\bf153}(2):1477-1486.
On-line at http:// www.springerlink.com/content/106500/}}

\author{Sh.~R.~Shakirov\\
Institute for Theoretical and Experimental Physics\\
Moscow, Russia\\
Moscow Institute of Physics and Technology\\
Dolgoprudny, Moscow Oblast, Russia}

\maketitle

\begin{abstract}
We propose a method for constructing systems of polynomial equations that
define submanifolds of degenerate binary forms of an arbitrary degeneracy
degree.
\end{abstract}


\section{Introduction}
\label{sec1}

This paper is devoted to studying degeneration of symmetric forms. As is known,
in the manifold of all symmetric forms on $\C^n$, degenerate forms constitute
a submanifold that is invariant under the action of the structure group
$SL(n)$. This submanifold can be defined by a single equation $D=0$, where
$D$, called the discriminant, is an $SL(n)$-invariant polynomial in the form
coefficients.

The discriminant submanifold has an interesting internal geometry: forms with
different degeneracy degrees produce a set of submanifolds in it, one
embedded into another. A natural problem is to construct systems of
polynomial equations that define these submanifolds. Such systems of
equations allow determining not only whether a given form is degenerate but
also how strongly degenerate it is.

The most general approach to solving this problem is described in~\cite{1} in
terms of representations of the group $SL(n)$. The simplest case $n=2$ is
related to the theory of algebraic equations in one variable. The present
paper is devoted to analyzing this case in detail.

Homogeneous symmetric forms and their degeneracy conditions often appear in
physical problems. In particular, a possible application is calculating
non-Gaussian generating functionals in field theory exactly
(nonperturbatively). In the formalism of continuum (functional) integration,
the integrals
$$
Z=\int e^{S_{ij}x^ix^j}\,d^nx=\const\cdot\frac1{\sqrt{\det S}\,},
$$
called Gaussian integrals, play an important role. They are well investigated
because $S$ is a quadratic form and its discriminant $D(S)$ is exactly equal
to its determinant $\det S$. The theory for calculating non-Gaussian
integrals (with forms of an arbitrary degree) is developed worse. It is very
likely related closely to the theory of discriminants and higher
discriminants.

\section{Binary forms and their roots}
\label{sec2}

A binary $k$-form is a symmetric form of degree $k$ on the space $\C^2$,
i.e., a homogeneous polynomial in two complex variables
$$
P_k(x,y)=a_kx^k+a_{k-1}x^{k-1}y+a_{k-2}x^{k-2}y^2+\dots+a_0y^k.
$$
The numbers $a_i$ are called the coefficients of the form $P$. Because the
fundamental field $\C$ is algebraically closed, any such polynomial can be
decomposed into a product of $k$ linear factors:
$$
P_k(x,y)=a_kx^k+a_{k-1}x^{k-1}y+\dots+a_0y^k=
(\alpha_1x+\beta_1y)(\alpha_2x+\beta_2y)\cdots
(\alpha_kx+\beta_ky).
$$
Consequently, the kernel of the form $P$ (the set of solutions of the
equation $P_k(x,y)=0$) consists of $k$ one-dimensional subspaces
$$
\Lambda_i=\bigl\{(x,y)\mid\alpha_ix+\beta_iy=0\bigr\}.
$$
These straight lines are called the roots of $P$. A binary $k$-form in
general position has exactly $k$ nonequal roots. Degenerate forms are those
with some roots equal to each other.

The roots of a binary form are one-dimensional subspaces in $\C^2$, i.e.,
points of the projective space $\C P^1$. A close relation between binary
$k$-forms and ordinary $k$th-degree algebraic equations consists in the
possibility to select some chart on $\C P^1$, for example, $y\ne0$. Because
$\dim\C P^1=1$, each root $\Lambda_i$ is represented in this chart by one
complex number $\lambda_i$. Moreover, the equation
$$
P_k(x,y)=a_kx^k+a_{k-1}x^{k-1}y+\dots+a_0y^k=0
$$
becomes
$$
P_k(z)=a_kz^k+a_{k-1}z^{k-1}+\dots+a_0=0,
$$
where $z=x/y$ is a coordinate on a projective space in the chart $y\ne0$.
This is a polynomial equation of degree $k$ in one complex variable, which
has exactly $k$ roots $\lambda_i$.

Geometrically, the space of all binary $k$-forms is a smooth manifold
$S^k\equiv S^k\C^2\sim\C^{k+1}$, in which degenerate forms constitute
projective submanifolds. The main one among them is the discriminant surface
$S^k_2$, the set of all $k$-forms with at least two equal roots. It is
defined by one polynomial equation
$$
D_k=0,
$$
where the polynomial $D_k$ is called the discriminant of $k$-forms. For
example, for a quadratic form
$$
P_2(x,y)=ax^2+bxy+cy^2,
$$
the discriminant is equal to
$$
D_2(a,b,c)=b^2-4ac,
$$
and for a cubic form
$$
P_3(x,y)=ax^3+bx^2y+cxy^2+dy^3,
$$
the discriminant is
$$
D_3(a,b,c,d)=-27a^2d^2+18abcd-4ac^3-4db^3+b^2c^2.
$$
Given the discriminant of a binary form, we can answer the question whether
the form is degenerate, i.e., whether it has a pair of coincident roots. Our
aim here is to construct higher discriminants, which would additionally allow
determining how many roots are equal and with what multiplicities.

On the discriminant surface $S^k_2$, there are submanifolds of forms with
higher degeneracy degrees. For example, $S^k_3$ is the set of all $k$-forms
with three equal roots, and $S^k_{2,2}$ is the set of all $k$-forms with two
pairs of equal roots.

For a given form, the degeneracy degree is determined by dividing the set of
its roots into coincident roots. For any partition
$$
k=k_1+k_2+\dots+k_p,
$$
where $k_1\ge k_2\ge\dots\ge k_p$, we define $S^k_{k_1,\dots,k_p}$ as the set
of $k$-forms with such type of coincident roots. To simplify calculations, we
omit the index $k_j$ in $S_{k_1,\dots,k_n}$ if it is equal to 1. For example,
the partition
$$
k=1+1+\dots+1
$$
corresponds to the space of all $k$-forms $S^k$, the partition
$$
k=2+1+\dots+1
$$
corresponds to the space of all degenerate $k$-forms $S^k_2$, and so on.

We must note that a partition $k=k_1+k_2+\dots+k_p$ does not imply that the
groups of coincident roots are distinct. For example, $S^k_4$ is a subset of
$S^k_{2,2}$ because the coincidence of four roots is a particular case of the
coincidence of two pairs of roots. In exactly the same way, $S^k_4$ is a
subset of $S^k_3$, and all these submanifolds are contained in $S^k_2$.

The submanifold $S^k_2\subset S^k$ is defined by one polynomial equation, but
submanifolds $S_{k_1,\dots,k_p}^k$ with $k_1>2$ have a lower dimension, and
should therefore be defined by systems of several polynomial equations. The
problem of describing such systems of equations was posed by Cayley as early
as the 19th century, but it still attracts interest and is investigated using
the most modern methods~\cite{2}.

The group of linear transformations $SL(2,\C)$
\begin{align*}
&x\to G_{11}x+G_{12}y,
\\[3mm]
&y\to G_{21}x+G_{22}y,
\end{align*}
where $G_{11}G_{22}-G_{12}G_{21}=1$, acts naturally on symmetric forms on
$\C^2$. Clearly, any submanifold of degenerate forms $S_{k_1,\dots,k_p}^k$,
just like the discriminant submanifold $S^k_2$, is invariant under this
action because the coincidence of a few lines on a plane is independent of
the choice of basis in this plane. Consequently, the discriminant of a binary
form is an $SL(2)$-invariant polynomial of its coefficients.

The $SL$ invariance is an important property of submanifolds of degenerate
forms and allows using representation theory to investigate them. This is
important because the kernel of a symmetric form is a discrete set of lines
only in the particular case of $\C^2$. For $n>2$, an arbitrary symmetric form
no longer decomposes into factors, in other words, there is no discrete set
of roots. It becomes impossible to formulate degeneration in terms of
coincidence of roots. A general theory of degeneration of symmetric forms
should most likely be formulated in group theory terms.

But in the case of binary forms, the kernel is equal to a discrete set of
roots, which can coincide with each other with different multiplicities. In
this case, we can construct equations defining forms with higher degeneracy
degrees using only the notion of roots, the Vieta theorem, and mathematical
logic.

\section{Discriminant}
\label{sec3}

From now on, a binary form $P(x,y)$ is considered an inhomogeneous polynomial
$P(z)$ in one complex variable $z$. The roots of the form correspond to the
roots of the algebraic equation $P(z)=0$.

We recall the main facts concerning algebraic equations of the $k$th degree
in one complex variable. Such an equation has the form
$$
a_kz^k+a_{k-1}z^{k-1}+\dots+a_0=0
$$
and has exactly $k$ solutions $\lambda_1,\lambda_2,\dots,\lambda_k$, which
are called its roots. The roots are related to the coefficients via the Vieta
formulas
$$
a_{k-i}=a_k(-1)^i\cdot\sigma_i(\lambda_1,\lambda_2,\dots,\lambda_k),
$$
where $\sigma_i$ denotes the $i$th elementary symmetric polynomial of roots,
\begin{align*}
\sigma_1(\lambda_1,\lambda_2,\dots,\lambda_k)&=
\lambda_1+\lambda_2+\dots+\lambda_k,
\\[3mm]
\sigma_2(\lambda_1,\lambda_2,\dots,\lambda_k)&=
\lambda_1\lambda_2+\lambda_1\lambda_3+\dots+\lambda_{k-1}\lambda_k,
\\[3mm]
&\;\,\vdots
\\[3mm]
\sigma_k(\lambda_1,\lambda_2,\dots,\lambda_k)&=
\lambda_1\lambda_2\cdots\lambda_k.
\end{align*}
The Vieta formulas follow at once from the decomposition
$$
a_kz^k+a_{k-1}z^{k-1}+\dots+a_0=
a_k(z-\lambda_1)(z-\lambda_2)\cdots(z-\lambda_k).
$$
Any symmetric polynomial of roots, as is known, can be represented as a
polynomial of elementary symmetric polynomials $\sigma_i$. Consequently, any
symmetric polynomial of the roots is actually a function of the coefficients.

A known method for finding the discriminant of a $k$th degree equation is to
consider the function of the roots
$$
D(\lambda_1,\lambda_2,\dots,\lambda_k)=\prod_{i<j}(\lambda_i-\lambda_j)^2,
$$
which is a symmetric polynomial of the roots and therefore a function of the
coefficients. Moreover, this polynomial is equal to zero if (and only if)
some two roots coincide. These properties correspond to the definition of
the discriminant. We conclude that the function $D$ should indeed coincide
with the discriminant up to a nonzero factor.

For example, for the quadratic equation $ax^2+bx+c=0$,
$$
D=(\lambda_1-\lambda_2)^2=\lambda_1^2-2\lambda_1\lambda_2+\lambda_2^2.
$$
Using the Vieta formulas $b/a=-(\lambda_1+\lambda_2)$ and
$c/a=\lambda_1\lambda_2$, we obtain
$$
D=\lambda_1^2+2\lambda_1\lambda_2+\lambda_2^2-4\lambda_1\lambda_2=
\biggl(\frac ba\biggr)^2-\frac{4c}{a}=\frac{b^2-4ac}{a^2}.
$$

Such a method for constructing the discriminant can be generalized to an
arbitrary type of root coincidence. It is necessary to write an expression in
the roots that vanishes {\it if and only if} the roots coincide in the
desired way (for example, three roots coincide). Here, it is convenient to
use the apparatus of mathematical logic.

\section{Higher discriminants}
\label{sec4}

We consider a $k$-form, its roots $\lambda_1,\lambda_2,\dots,\lambda_k$, and
the logical algebra $\Lambda$ generated by all logical statements
$$
E_{ij}:(\lambda_i=\lambda_j)
$$
together with the logical operations AND and OR, respectively denoted by $+$
and $\times$.

For any type of coincidence $S_{k_1,\dots,k_p}^k$, we let
$L[S_{k_1,\dots,k_p}^k]$ denote the logical expression built from elementary
expressions $E_{ij}$, which is called the definition of this type of
coincidence. The definition is constructed very naturally for any type of
coincidences. For example, the definition of the discriminant submanifold
$S^k_2$ is that some pair of roots coincides: either the first root is equal
to the second, or the first root is equal to the third, and so on. We thus
obtain
$$
L[S^k_2]=E_{12}\times E_{13}\times\dots\times E_{k-1,k}=
\underset{i<j}{\times} E_{ij}.
$$
Some simple examples of definitions are
\begin{align*}
&L[S^3_3]=E_{12}+E_{13}+E_{23},
\\[3mm]
&L[S^4_{2,2}]=(E_{12}+E_{34})\times(E_{13}+E_{24})\times(E_{14}+E_{23}),
\\[3mm]
&L[S^4_4]=E_{12}+E_{13}+E_{14}+E_{23}+E_{24}+E_{34},
\\[3mm]
&L[S^4_3]=(E_{12}+E_{13}+E_{23})\times(E_{12}+E_{14}+E_{24})\times{}
\\[3mm]
&\hphantom{L[S^4_3]={}}\times(E_{23}+E_{24}+E_{34})\times
(E_{13}+E_{14}+E_{34}).
\end{align*}

There is a direct relation between the logical expressions and equations: for
each statement in $\Lambda$, there exists an equivalent system of equations
imposed on roots. For example, an elementary statement $E_{ij}$ corresponds
to the equation (squared for symmetry)
$$
(\lambda_i-\lambda_j)^2=0.
$$
The operation OR corresponds to multiplying the operands because a product is
equal to zero if one of the factors is zero. The operation AND literally
corresponds to taking a system of equations. For example,
\begin{align*}
&E_{12}\mapsto(\lambda_1-\lambda_2)^2=0,
\\[3mm]
&E_{34}\mapsto(\lambda_3-\lambda_4)^2=0,
\\[3mm]
&E_{12}\times E_{34}\mapsto
(\lambda_1-\lambda_2)^2(\lambda_3-\lambda_4)^2=0,
\\[3mm]
&E_{12}+E_{34}\mapsto\begin{cases}
(\lambda_1-\lambda_2)^2=0,\\[2mm]
(\lambda_3-\lambda_4)^2=0.\end{cases}
\end{align*}
We note that the discriminant case $S^k_2$ corresponds to one equation,
$$
L[S^k_2]=E_{12}\times E_{13}\times\dots\times E_{k-1,k}\mapsto
\prod_{i<j}(\lambda_i-\lambda_j)^2=0,
$$
which was already mentioned above. The coincidence of three roots for a
3-form corresponds to the system of equations
$$
L[S^3_3]=E_{12}+E_{13}+E_{23}\mapsto\begin{cases}
(\lambda_1-\lambda_2)^2=0,\\[2mm]
(\lambda_1-\lambda_3)^2=0,\\[2mm]
(\lambda_2-\lambda_3)^2=0.\end{cases}
$$

At first glance, the logical expression
$$
L[S^4_{2,2}]=(E_{12}+E_{34})\times(E_{13}+E_{24})\times(E_{14}+E_{23})
$$
has no corresponding system of equations. In this and similar cases, we must
expand the expression using the ordinary distributive law for $\times$ and
$+$,
$$
a\times(b+c)=a\times b+a\times c.
$$
The definition of $L[S^4_{2,2}]$ then becomes
\begin{align*}
L[S^4_{2,2}]={}&E_{12}\times E_{13}\times E_{14}+
E_{12}\times E_{24}\times E_{14}+{}
\\[3mm]
&{}+E_{12}\times E_{13}\times E_{23}+
E_{12}\times E_{24}\times E_{23}+E_{34}\times E_{13}\times E_{14}+{}
\\[3mm]
&{}+E_{34}\times E_{24}\times E_{14}+E_{34}\times E_{13}\times E_{23}+
E_{34}\times E_{24}\times E_{23}.
\end{align*}
The corresponding system of equations is
\begin{align*}
&(\lambda_1-\lambda_2)^2(\lambda_1-\lambda_3)^2(\lambda_1-\lambda_4)^2=0,
\\[3mm]
&(\lambda_1-\lambda_2)^2(\lambda_2-\lambda_4)^2(\lambda_1-\lambda_4)^2=0,
\\[3mm]
&(\lambda_1-\lambda_2)^2(\lambda_1-\lambda_3)^2(\lambda_2-\lambda_3)^2=0,
\\[3mm]
&(\lambda_1-\lambda_2)^2(\lambda_2-\lambda_4)^2(\lambda_2-\lambda_3)^2=0,
\\[3mm]
&(\lambda_3-\lambda_4)^2(\lambda_1-\lambda_3)^2(\lambda_1-\lambda_4)^2=0,
\\[3mm]
&(\lambda_3-\lambda_4)^2(\lambda_2-\lambda_4)^2(\lambda_1-\lambda_4)^2=0,
\\[3mm]
&(\lambda_3-\lambda_4)^2(\lambda_1-\lambda_3)^2(\lambda_2-\lambda_3)^2=0,
\\[3mm]
&(\lambda_3-\lambda_4)^2(\lambda_2-\lambda_4)^2(\lambda_2-\lambda_3)^2=0.
\end{align*}

For any type of coincidences, the method described above allows constructing
a system of equations imposed on the roots that is equivalent to the
definition of this type of coincidence, and it consequently defines the space
of forms with that type of coincidence.

The final step is to express this system of equations in terms of the
coefficients of the form. We suppose that the system of equations obtained
for some type of coincidence has the form
\begin{align*}
P_1(\lambda_1,\lambda_2,\dots,\lambda_k)&=0,
\\[3mm]
&\;\,\vdots
\\[3mm]
P_m(\lambda_1,\lambda_2,\dots,\lambda_k)&=0.
\end{align*}
Generally, the polynomials $P_i$ are not symmetric functions of the roots,
and they therefore cannot be immediately expressed as functions of the
coefficients. This problem is easily solved by symmetrizing this system,
i.e., by passing to the equivalent system
\begin{align*}
&P_1+P_2+\dots+P_m=0,
\\[3mm]
&P_1P_2+P_1P_3+\dots+P_{m-1}P_m=0,
\\[3mm]
&\quad\vdots
\\[3mm]
&P_1P_2\dots P_m=0.
\end{align*}
The equations in this system are symmetric in the polynomials $P_i$.
Consequently, they are symmetric in the roots and can be expressed in terms
of the coefficients via the Vieta formulas. As a result, we have a solution
of the problem: a system of polynomial equations for the coefficients of the
form is obtained, defining a submanifold of forms with the given type of
degeneracy.

\section{The case $S^3_3$}
\label{sec5}

To illustrate the method, we consider conditions for the coincidence of three
roots of a cubic binary form
$$
P_3(x,y)=ax^3+bx^2y+cxy^2+dy^3
$$
and the corresponding set of degenerate forms $S^3_3$. The definition in this
case is that all three roots are equal and
$$
L[S^3_3]=E_{12}+E_{13}+E_{23}.
$$
The corresponding system of equations
\begin{align*}
&(\lambda_1-\lambda_2)^2=0,
\\[3mm]
&(\lambda_1-\lambda_3)^2=0,
\\[3mm]
&(\lambda_2-\lambda_3)^2=0
\end{align*}
after symmetrization becomes
\begin{align*}
&(\lambda_1-\lambda_2)^2+(\lambda_1-\lambda_3)^2+(\lambda_2-\lambda_3)^2=0,
\\[3mm]
&(\lambda_1-\lambda_3)^2(\lambda_1-\lambda_2)^2+
(\lambda_1-\lambda_3)^2(\lambda_2-\lambda_3)^2+
(\lambda_1-\lambda_2)^2(\lambda_2-\lambda_3)^2=0,
\\[3mm]
&(\lambda_1-\lambda_2)^2(\lambda_1-\lambda_3)^2(\lambda_2-\lambda_3)^2=0
\end{align*}
and can be expressed in terms of $a$, $b$, $c$, and $d$ using the Vieta
formulas
\begin{align*}
&\frac{b}{a}=-(\lambda_1+\lambda_2+\lambda_3),
\\[3mm]
&\frac{c}{a}=\lambda_1\lambda_2+\lambda_1\lambda_3+\lambda_2\lambda_3,
\\[3mm]
&\frac{d}{a}=-\lambda_1\lambda_2\lambda_3.
\end{align*}
We omit the simple calculations leading to the result
\begin{align*}
&2(b^2-3ac)=0,
\\[3mm]
&(b^2-3ac)^2=0,
\\[3mm]
&-27a^2d^2+18abcd-4ac^3-4db^3+b^2c^2=0.
\end{align*}
This system of equations is equivalent to the definition of the coincidence
of three roots, and it therefore defines the submanifold $S^3_3\subset S^3$.
The second equation is dependent on the first and can be omitted.

We must note that this system is not $SL(2)$-invariant, although the
submanifold $S^3_3$ defined by this system is invariant. This means that
another system, transformed under the $SL(2)$ action, defines the same
manifold. If we let $G$ denote an arbitrary element of $SL(2)$, then we
obtain infinitely many systems
\begin{align*}
&2Z(G)=0,
\\[3mm]
&\bigl(Z(G)\bigr)^2=0,
\\[3mm]
&-27a^2d^2+18abcd-4ac^3-4db^3+b^2c^2=\inv=0,
\end{align*}
where $Z(G)$ is the orbit of $b^2-3ac$ under the group action,
$$
Z(G)=G_{11}^2(b^2-3ac)+G_{11}G_{12}(bc-9ad)+G_{12}^2(c^2-3bd).
$$
Any of these systems define the submanifold $S^3_3$ of forms with three
coincident roots.

An interesting fact is that the set of common zeros of the whole orbit
$Z(G)$, i.e., the set of forms for which the whole orbit $Z(G)$ vanishes, is
exactly $S^3_3$. This allows regarding $Z(G)$ as a higher discriminant for
the case of the coincidence of three roots out of three. The whole orbit, of
course, is an $SL(2)$-invariant object. From the representation theory
standpoint, this orbit realizes an irreducible representation of the $SL(2)$
group in the space of quadratic polynomials of coefficients.

\section{The case $S^4_3$}
\label{sec6}

As a second example, we construct a system of equations defining forms of
degree four,
$$
P_4(x,y)=ax^4+bx^3y+cx^2y^2+dxy^3+ey^4,
$$
which have one root of multiplicity three. The definition of the considered
coincidence is
\begin{align*}
L[S^4_3]={}&(E_{12}+E_{13}+E_{23})\times(E_{12}+E_{14}+E_{24})\times{}
\\
&{}\times(E_{23}+E_{24}+E_{34})\times(E_{13}+E_{14}+E_{34}).
\end{align*}
Using standard rules of mathematical logic, we can easily verify that this
definition is equivalent to the simpler logical statement
$$
L'=E_{12}\times E_{34}+E_{13}\times E_{24}+E_{14}\times E_{23},
$$
which corresponds to the system of equations
\begin{align*}
&(\lambda_1-\lambda_2)^2(\lambda_3-\lambda_4)^2=0,
\\[3mm]
&(\lambda_1-\lambda_3)^2(\lambda_2-\lambda_4)^2=0,
\\[3mm]
&(\lambda_2-\lambda_3)^2(\lambda_1-\lambda_4)^2=0.
\end{align*}
Symmetrizing this system, we obtain
\begin{align*}
&(\lambda_1-\lambda_2)^2(\lambda_3-\lambda_4)^2+
(\lambda_1-\lambda_3)^2(\lambda_2-\lambda_4)^2+
(\lambda_2-\lambda_3)^2(\lambda_1-\lambda_4)^2=0,
\\[3mm]
&(\lambda_1-\lambda_2)^2(\lambda_3-\lambda_4)^2
(\lambda_1-\lambda_3)^2(\lambda_2-\lambda_4)^2+
(\lambda_1-\lambda_2)^2(\lambda_3-\lambda_4)^2\times{}
\\[3mm]
&\qquad\times
(\lambda_2-\lambda_3)^2(\lambda_1-\lambda_4)^2+
(\lambda_1-\lambda_3)^2(\lambda_2-\lambda_4)^2
(\lambda_2-\lambda_3)^2(\lambda_1-\lambda_4)^2=0,
\\[3mm]
&(\lambda_1-\lambda_2)^2(\lambda_3-\lambda_4)^2
(\lambda_1-\lambda_3)^2(\lambda_2-\lambda_4)^2
(\lambda_2-\lambda_3)^2(\lambda_1-\lambda_4)^2=0.
\end{align*}
This system can be expressed in terms of the coefficients $a$, $b$, $c$, $d$,
and $e$ using the Vieta formulas
\begin{align*}
&\frac{b}{a}=-(\lambda_1+\lambda_2+\lambda_3+\lambda_4),
\\[3mm]
&\frac{c}{a}=\lambda_1\lambda_2+\lambda_1\lambda_3+
\lambda_1\lambda_4+\lambda_2\lambda_3+\lambda_2\lambda_4+\lambda_3\lambda_4,
\\[3mm]
&\frac{d}{a}=-(\lambda_1\lambda_2\lambda_3+\lambda_1\lambda_2\lambda_4+
\lambda_1\lambda_3\lambda_4+\lambda_2\lambda_3\lambda_4),
\\[3mm]
&\frac{e}{a}=\lambda_1\lambda_2\lambda_3\lambda_4.
\end{align*}
This calculation leads to the result
\begin{align*}
&2(c^2-3bd+12ae)=0,
\\[3mm]
&(c^2-3bd+12ae)^2=0,
\\[3mm]
&D_4=0.
\end{align*}

The left-hand side of the third equation is equal to the discriminant of a
fourth-degree form. As in the preceding example, the second equation is
dependent on the first and can be omitted.

An interesting fact is that all equations in this system are invariant under
the $SL(2)$ action. We call the invariant polynomial $c^2-3bd+12ae$ the
apolara.
We find that the higher discriminant that is relevant for this case of
the coincidence of three roots out of four is a system of two invariants:
the discriminant and the apolara.

\section{Conclusion}
\label{sec7}

The proposed method allows obtaining systems of equations that define binary
forms with an arbitrary partition of the roots into coinciding roots. The
idea of the method is to select a chart on the projective space $\C P^1$
where roots are represented by complex numbers and the form itself is
represented by a polynomial in one variable. The required system of equations
is then constructed in the space of roots and finally expressed in terms of
the coefficients using the symmetric Vieta formulas.

The questions of the correspondence between the orbits of the group $SL(2)$
and the conditions for the degeneracy of the forms and also of the
generalization of this correspondence to the case of space dimensions $n>2$
are still unresolved.

\subsection*{Acknowledgments}
The author is grateful to V.~Dolotin and A.~Morozov for the very useful
discussions.

\smallskip

This work was supported in part by the Federal Nuclear Energy Agency, the
Program for Supporting Leading Scientific Schools
(Grant No.~NSh-8004.2006. 2), and the Russian Foundation for Basic Research
(Grant Nos.~RFBR--Italy 06-01-92059-CE and~07-02-006454).

\end{document}